\documentclass{birkjour}
\usepackage{amsfonts}
\usepackage{amsmath}
\usepackage{amssymb}
\usepackage{mathrsfs} 
\usepackage{color} 
\usepackage[english]{babel}
\usepackage[T1]{fontenc}
\usepackage{amsthm}
\input xy 
\xyoption{all}
\usepackage{enumitem}
\usepackage{hyperref}
\hypersetup{colorlinks=true}
\usepackage{textcomp}

\theoremstyle{plain}
\newtheorem{thm}{Theorem}
\newtheorem*{thm*}{Theorem}
\theoremstyle{definition}
\newtheorem{lem}[thm]{Lemma}
\newtheorem{propo}[thm]{Proposition}
\newtheorem{rmk}[thm]{Remark}

\newcommand{\PU}{\gr{P}(U)}
\newcommand{\PV}{\gr{P}(V)}

\newcommand{\mappa}[2]{\xymatrix@1{{#1} \ar[r] & {#2}}}
\newcommand{\spazio}{\rule{1 pt}{0 cm}}
\newcommand{\gr}[1]{\mathbf{#1}}

\newcommand{\OPU}{\mathcal{O}_{\PU}}

\newcommand{\GRA}{\gra(m,\Lambda^{2}V)}
\newcommand{\GRAd}{\gra(2,\Lambda^{2}V)}

\newcommand{\PLVs}{\gr{P}(\Lambda^{2}V^{*})}
\newcommand{\PLV}{\gr{P}(\Lambda^{2}V)}
\newcommand{\GG}{\mathbb{G}}
\newcommand{\Pnm}{\gr{P}(V)}

\renewcommand{\Bbb}{\mathbb}
\renewcommand{\bold}{\mathbf}

\DeclareMathOperator{\Deg}{Deg}

\DeclareMathOperator{\Pf}{Pf}
\DeclareMathOperator{\PGL}{PGL}
\DeclareMathOperator{\GL}{GL}

\DeclareMathOperator{\gra}{\textbf{Gr}}
\DeclareMathOperator{\grapr}{\GG r}

\DeclareMathOperator{\im}{Im}

\DeclareMathOperator{\Vi}{V}
\DeclareMathOperator{\Hh}{H}
\DeclareMathOperator{\hh}{h}

\begin{document}

\title[Degeneracy loci of twisted differential forms and linear line complexes]{Degeneracy loci of twisted differential forms\\ and linear line complexes}
\author{Fabio Tanturri}
\address{Mathematik und Informatik\\
Geb\"aude E.2.4\\
Universit\"at des Saarlandes\\
D-66123 Saarbr\"ucken\\
Germany}
\email{tanturri@math.uni-sb.de}
\thanks{Partially supported by the PRIN 2010/2011 ``Geometria delle variet\`a algebriche''}



\date{\today}

\begin{abstract}
We prove that the Hilbert scheme of degeneracy loci of pairs of global sections of $\Omega_{\mathbb{P}^{n-1}}^{}(2)$, the twisted cotangent bundle on $\mathbb{P}^{n-1}$, is unirational and dominated by the Grassmannian of lines in the projective space of skew-symmetric forms over a vector space of dimension $n$. We provide a constructive method to find the fibers of the dominant map. In classical terminology, this amounts to giving a method to realize all the pencils of linear line complexes having a prescribed set of centers. In particular, we show that the previous map is birational when $n=4$.
\end{abstract}

\maketitle

%%%%%%%%%%%%%%%%%%%%%%%%%%%%%%%%%%%%%%%%%%%%%%%%%%%%%%%%%%%%%%%%%%%%%%%%%%%%%%%%%%%%%%%%%%%%%%%%%%%
%%%%%%%%%%%%%%%%%%%%%%%%%%%%%%%%%%%%%%%%%%%%%%%%%%%%%%%%%%%%%%%%%%%%%%%%%%%%%%%%%%%%%%%%%%%%%%%%%%%
\section{Introduction}
Degeneracy loci of morphisms of the form $\phi:\mappa{\mathcal{O}_{\mathbb{P}^{n-1}}^{m}}
{\Omega_{\mathbb{P}^{n-1}}^{}(2)}$ arise naturally in algebraic geometry 
and have been extensively studied. In classical terminology, a general degeneracy locus of this form is the set of centers of complexes belonging to a general linear system of dimension $m-1$ of linear line complexes in $\mathbb{P}^{n-1}$, i.e.~hyperplane sections of the Grassmannian $\grapr(1,\mathbb{P}^{n-1})$ of lines in $\mathbb{P}^{n-1}=\PV$ embedded in $\PLV$ via the Pl\"ucker map (cfr.~Sect.~\ref{linlincom}). This nice geometric interpretation, together with the fact that many classical algebraic varieties arise this way, motivated the interest of many classical algebraic geometers as Castelnuovo, Fano, and Palatini. Nonetheless, these loci have been recently studied by several authors, for example Chang and Ottaviani. A more detailed historical account and a few classical examples can be found, for instance, in \cite{BazanMezzetti,FaenziFania}.

In general, in order to parameterize the possible degeneracy loci $X_{\phi}$ as $\phi$ varies, it is useful to take a modern approach and introduce $\mathcal{H}$ as the union of the irreducible 
components, in the Hilbert scheme of subschemes of $\mathbb{P}^{n-1}$, containing $X_{\phi}$ for general $\phi$'s.

Let $\PV$ be the projectivization of an $n$-dimensional vector space $V$. By the Euler sequence we can
identify a morphism of the form above with a skew-symmetric matrix of linear
forms in $m$ variables, or with an $m$-tuple of elements in $\Lambda^{2}V$; the natural 
$\GL_{m}$-action does not modify its degeneracy locus, so we 
get the rational map
\begin{equation}
\label{rhointro}
\rho:
\xymatrix{
\GRA \ar@{-->}[r] & \mathcal{H}
}
\end{equation}
sending $\phi$ to $X_{\phi}$.

The behavior of the map $\rho$ is fully understood for most of the values of $(m,n)$. In \cite{TanturriDegeneracy} it is shown that it is birational whenever $4 \leq m < n-1$ or $(m,n)=(3,5)$; if $(m,n)=(3,6)$, $\rho$ is dominant but $4:1$, while for $m=3$ and $n>6$ it is generically injective but not dominant. Partial results in this line of thought were already provided in \cite{BazanMezzetti,FaniaMezzetti,FaenziFania}.

\smallskip We are interested here in the case $m=2$. Our aim is to give a description of the behavior of the map $\rho$, as $n$ varies. The case $n=6$ has already been tackled in \cite{BazanMezzetti}, where $\rho$ was proved to be dominant and its fibers were described and shown to be two-dimensional. Our main result is the following

\begin{thm*}
Let $n \in \mathbb{N}$ such that $n \geq 4$, $V$ be a vector space of dimension $n$ and let 
\[\rho:
\xymatrix{
\GRAd \ar@{-->}[r] & \mathcal{H} 
}\]
be the rational morphism introduced in \eqref{rhointro}, sending the class of a morphism $\phi:\mappa{\mathcal{O}_{\PV}^{2}}
{\Omega^{}_{\PV}(2)}$ to its degeneracy locus $X_{\phi}$, considered as a point in the Hilbert scheme. Then $\rho$ is dominant. Moreover
\begin{enumerate}[label=\textup{\roman{*}.}, ref=(\roman{*})]
\item if $n$ is even, the general element of $\mathcal{H}$ is the union of $\frac{n}{2}$ lines spanning the whole $\PV$. The general fiber is the Grassmannian $\grapr(1,\sigma)$ of lines lying on a suitable $(\frac{n-2}{2})$-dimensional projective space $\sigma$. In particular, $\rho$ is birational for $n=4$;
\item if $n$ is odd, the general element of $\mathcal{H}$ is the image in $\PV$ of a map
\[
\xymatrix{
\mathbb{P}^{1} \ar[rr]^-{[f_{1}:\dotso:f_{n}]} && \PV,
}
\]
where $f_{1},\dotsc,f_{n}$ are forms of degree $\frac{n-1}{2}$ spanning the whole vector space $\Hh^{0}(\mathbb{P}^{1},\mathcal{O}_{\mathbb{P}^{1}}(\frac{n-1}{2}))$. The general fiber of $\rho$ has dimension $\frac{n^{2}-3n}{2}$.
\end{enumerate}
\end{thm*}

When $n$ is even, the general $X_{\phi}$ is union of $\frac{n}{2}$ lines. A crucial step in proving the first part of the statement is the reinterpretation of a morphism $\phi$ as a pencil $\ell_{\phi}$ of linear line complexes. The $\frac{n}{2}$ lines are the set of centers of the linear line complexes belonging to $\ell_{\phi}$; using a similar approach to the one adopted in \cite{BazanMezzetti}, we are able to generalize the existing result to all even values of $n$.

When $n$ is odd, the general $X_{\phi}$ is a rational curve, image of the map $[f_{1}:\dotso:f_{n}]$ as in the statement. A priori, the forms $f_{1},\dotsc,f_{n}$ have to be the Pfaffians of order $n-1$ of a linear $n \times n$ skew-symmetric matrix, but we will show that a general $n$-tuple of forms of degree $\frac{n-1}{2}$ arises as such.

As pointed out later, an explicit construction of the preimages of a given degeneracy locus is possible. For the special case $n=4$, this amounts to the following: it is possible to construct the unique pencil of linear line complexes having two given skew lines as set of centers.

These results complete the general picture, describing the behavior of the map $\rho$ for arbitrary values of $m$. For $m>3$, such description had already been given in \cite{TanturriDegeneracy}; the above theorem, together with the previous contributions, shows that the unique pairs $(m,n)$ with $2 \leq m < n-1$ for which $\rho$ is not generically injective are $(2,n)$ with $n\geq 5$ and $(3,6)$, while the unique pairs for which $\rho$ is not dominant are $(3,n)$ with $n\geq 7$.

We remark that the results contained in this paper hold over non-necessarily algebraically closed fields of arbitrary characteristic.
%%%%%%%%%%%%%%%%%%%%%%%%%%%%%%%%%%%%%%%%%%%%%%%%%%%%%%%%%%%%%%%%%%%%%%%%%%%%%%%%%%%%%%%%%%%%%%%%%%%
%%%%%%%%%%%%%%%%%%%%%%%%%%%%%%%%%%%%%%%%%%%%%%%%%%%%%%%%%%%%%%%%%%%%%%%%%%%%%%%%%%%%%%%%%%%%%%%%%%%
\section{Notation and preliminaries}

\subsection{Notation and geometric interpretation}
Let $\mathbf{k}$ be any field and let $n \in \mathbb{N}$ such that $n \geq 4$. We will denote by $U, V$ two $\mathbf{k}$-vector spaces of dimensions $2$, $n$ respectively; by $\PU$ and $\PV$ we will mean the projective spaces of their 1-quotients, i.e.~$\Hh^{0}(\PU,\OPU(1))\cong U$. We set $\{y_{0},y_{1}\}$ and $\{x_{0},\dotsc,x_{n-1}\}$ to be the bases of $U$ and $V$ respectively.

Let $\Omega_{\PV}$ be the rank $n-1$ vector bundle of differential forms on $\PV$ and let $\phi:U \otimes \mathcal{O}_{\PV}\rightarrow \Omega_{\PV}(2)$ be a general morphism. We can define $X=X_{\phi}$ to be the degeneracy locus associated to $\phi$, i.e.~the scheme cut out by the maximal minors of the matrix locally representing $\phi$.

By the Euler sequence, the map $\phi$ can be interpreted also as an $(n \times n)$ skew-symmetric matrix $N_{\phi}=y_{0}N_{0}+y_{1}N_{1}$ of linear forms in $y_{0},y_{1}$. A point $p$ in $V$ represents a point in $X_{\phi}$ if and only if $p \in \ker(b_{0}N_{0}+b_{1}N_{1})$ for some $(b_{0},b_{1})\neq (0,0)$ (cfr.~\cite[\textsection 3.2]{Ottaviani}). Thus, the geometry of $X_{\phi}$ strongly depends on the parity of $n$: when $n$ is even, it is a scroll over the hypersurface cut out by the Pfaffian of $N_{\phi}$, i.e.~a union of $\frac{n}{2}$ lines in $\PV$; when $n$ is odd, it is the blow-up of $\mathbb{P}^{1}$ along the set of points described by the $n$ Pfaffians of order $n-1$ of $N_{\phi}$, which is empty for a general $\phi$. For more details about this geometric interpretation, we refer to \cite[\textsection 2]{TanturriDegeneracy}.

\subsection{Linear line complexes} \label{linlincom}
\noindent Another interpretation (cfr.~\cite{BazanMezzetti}) of a morphism $\phi:\mappa{\mathcal{O}_{\PV}^{2}}{\Omega_{\PV}(2)}$ is the following: let $\GG:=\grapr(1,\PV)$ be the Grassmannian of lines in $\PV$, 
embedded in $\PLV$ via the Pl\"ucker map. The dual space $\PLVs$ parameterizes hyperplane sections of $\GG$,
or, in classical terminology, \emph{linear line complexes} in $\PV$.

An element $A \in \Lambda^{2} V$, up to constants, may be regarded as an element of $\PLVs$, hence giving rise to a linear line complex $\Gamma$. A point $p \in \PV$
is called a \emph{center} of $\Gamma$ if all the lines through $p$ belong to $\Gamma$; the \emph{singular space} of $\Gamma$, i.e.~$\Deg(A):=\mathbf{P}(\ker(A))$, turns out to be the set of centers of $\Gamma$.

Let $\ell \in \GG$ and $\mathbb{T}_{\ell} \GG$ the corresponding tangent space. The hyperplane $\Vi(A)$ in $\PLV$ contains
$\mathbb{T}_{\ell} \GG$, i.e.~$A$ belongs to the dual variety $\check{\GG}$, if and only if $\Deg(A) \supseteq \ell$.

We distinguish the following two cases:
\begin{itemize}
\item if $n$ is even, then a general linear line complex $\Gamma$ has no center, as a general
skew-symmetric matrix of even order has maximal rank. Linear line complexes $\Gamma$ corresponding
to the points of $\check{\GG}$ have at least a line as center and will be called \emph{special}; 
\item if $n$ is odd, then a general linear line complex $\Gamma$ has a point as center. $\Gamma$ will be said to be \emph{special} if it corresponds to a point of $\check{\GG}$: in this case, the center of
$\Gamma$ contains (at least) a $\mathbb{P}^{2}$.
\end{itemize}

We can therefore interpret the degeneracy locus of a general morphism $\phi:\mappa{\mathcal{O}_{\PV}^{2}}{\Omega_{\PV}(2)}$ as
\[
\Deg(N_{\phi}):=\bigcup_{A \in N_{\phi}} \Deg(A), 
\]
where the skew-symmetric matrix $N_{\phi}=y_{0}N_{0}+y_{1}N_{1}$ is regarded as a line in $\PLVs$. 
Special complexes on $N_{\phi}$ are parameterized by the intersection $N_{\phi} \cap \check{\GG}$.

%%%%%%%%%%%%%%%%%%%%

%%%%%%%%%%%%%%%%%%%%%%%%%%%%%%%%%%%%%%%%%%%%%%%%%%%%%%%%%%%%%%%%%%%%%%%%%%%%%%%%%%%%%%%%%%%%%%%%%%%
%%%%%%%%%%%%%%%%%%%%%%%%%%%%%%%%%%%%%%%%%%%%%%%%%%%%%%%%%%%%%%%%%%%%%%%%%%%%%%%%%%%%%%%%%%%%%%%%%%%
\section{\texorpdfstring{The behavior of $\rho$: even case}{The behavior of rho: even case}}
Within this section, we assume $n$ even. We will show that the map $\rho$
\[\rho:
\xymatrix{
\GRAd \ar@{-->}[r] & \mathcal{H} 
}\]
defined in \eqref{rhointro} is dominant. This leads naturally to asking what is the preimage of a general point $X_{\phi}$ in $\mathcal{H}$, i.e.~which lines $N \subset \PLVs$ have $\Deg(N)=X_{\phi}$. We will give a geometric description of the fibers, providing a constructive procedure to realize the elements of a general preimage.

Recall that a linear line complex $A \in \PLVs$ is special if its center contains a $\mathbb{P}^{1}$. 
We distinguish special complexes \emph{of the first type}, whose center is exactly a line, from 
special complexes \emph{of the second type}, whose center is at least a $\mathbb{P}^{3}$.

A general matrix $N=N_{\phi}$ with linear forms in $\mathbf{k}[y_{0},y_{1}]$ as entries
has corank two in $\frac{n}{2}$ distinct points of $\mathbb{P}^{1}$, corresponding to the roots of $\Pf(N)$; $N$ has maximal rank in any other point. This means that the general line of complexes $N \subset \PLVs$ does not contain any special complex of the second type, and that $\Deg(N)$ is the union of $\frac{n}{2}$ lines $\{\ell_{1},\dotsc,\ell_{\frac{n}{2}}\}$.

We claim that these lines are general, in the sense that their span is the whole $\PV$.
It is clear that this condition is open, so it is sufficient to exhibit a matrix $N$ satisfying it. We will do more, showing constructively in Proposition \ref{bohbohboh} that any set of lines $\{\ell_{1},\dotsc,\ell_{\frac{n}{2}}\}$ spanning $\PV$ is $\Deg(N)$ for some pencil of complexes $N$.

Let us examine the Gauss map
\[
\zeta:\xymatrix{\check{\GG} \ar@{-->}[r] & \GG}
\]
which sends a special complex of the first type $A$ to the point in $\GG$ corresponding to the line $\Deg(A)$. Fixing a line $\ell \in \GG$, the fiber $\zeta^{-1}(\ell)$ is a linear space. A complex $A \in \GG$ has $\ell$ as center if and only if the hyperplane $\Vi(A)$ in $\PLV$ contains $\mathbb{T}_{\ell}\GG$, so the space of such $A$'s has dimension
\[
\dim(\zeta^{-1}(\ell))=
\dim(\PLV)-\dim(\GG)-1=
\frac{1}{2}(n-1)(n-4).
\]

We observe that, given a linear space $S$ of dimension $n-3$ in $\PV$, it is uniquely determined a complex $H$, up to constants, such that the center of $H$ is $S$. Indeed, all the lines in $\PV$ intersecting $S$ have to be contained in $\Gamma=\Vi(H)\cap \GG$, but this is a linear condition in the Pl\"ucker coordinates.

Fixing a set of lines $\{\ell_{1},\dotsc,\ell_{\frac{n}{2}}\}$ spanning $\PV$, for any $j$ such that $1 \leq j \leq \frac{n}{2}$ we will denote by $H_{j}\in \PLVs$ the unique complex having ${<}\ell_{i}{>}_{i\neq j}$ as center. We will denote by $F_{i}$ the $(\frac{n}{2}-2)$-dimensional linear span ${<}H_{j}{>}_{j\neq i}\subset\PLVs$; a point in $F_{i}$ is a complex having at least $\ell_{i}$ as center.

\begin{rmk}
\label{possoassumere}
Let $\ell_{1},\dotsc,\ell_{\frac{n}{2}}$ be lines spanning $\PV$. Up to a projective transformation of $\PV$, we may assume that
\[
\forall \, i, \quad \ell_{i}=\bigcap_{j \notin \{2i-2,2i-1\}} \Vi(x_{j}^{*}), %\cap \bigcap_{j > 2i-1} \Vi(x_{j}^{*})
\]
being $\Vi(x_{j}^{*})$ the hyperplane in $\PV$ whose points have $x_{j}$ coordinate zero.
%\begin{array}{c}
%\ell_{1}={<}[1:0:0:\dotso],[0:1:0:\dotso]{>}\\
%\ell_{2}={<}[0:0:1:0:0:\dotso],[0:0:0:1:0:\dotso]{>}\\
%\vdots \\
%\ell_{\frac{n}{2}}={<}[0:\dotso:0:1:0],[0:\dotso:0:0:1]{>}.
%\end{array}
%\]
With this choice, $H_{i}$ is the complex represented by a skew-symmetric matrix with $(j,k)$-th entry
\[
\left\{
\begin{array}{ll}
\alpha & \mbox{if } j+1=2i=k\\
-\alpha & \mbox{if } k+1=2i=j\\
0 & \mbox{otherwise}
\end{array}
\right.
\]
for some $\alpha \in \mathbf{k}\setminus\{0\}$. The elements of the fiber $\zeta^{-1}(\ell_{i})$ have zero entries in the {$(2i-1)$-th}, $(2i)$-th rows and columns.
\end{rmk}

\begin{propo}
\label{bohbohboh}
Let $L$ be the union of lines $\ell_{1},\dotsc,\ell_{\frac{n}{2}}$ spanning $\PV$. Let $\sigma$ be the $(\frac{n-2}{2})$-dimensional linear space ${{<}}H_{i}{{>}}_{1 \leq i \leq \frac{n}{2}} \subset \PLVs$. Then, for a line $N \subset \PLVs$, the following are equivalent:
\begin{enumerate}[label=\textup{\roman{*}.}, ref=(\roman{*})]
\item $N \subseteq \sigma$ and $N \cap F_{i} \cap F_{j}=\emptyset$ for any $1 \leq i < j \leq \frac{n}{2}$;
\item $N$ contains no special complexes of the second type and ${\Deg(N)=L}$.%\{\ell_{1},\dotsc,\ell_{\frac{n}{2}}\}$.
\end{enumerate}
\begin{proof} \spazio
\begin{itemize}
\item[i. $\Rightarrow$ ii.]
Let $N \subseteq \sigma$. For any $j$, the linear space $F_{j}$ has dimension $\frac{n}{2}-2$, so any line $N \subseteq \sigma$ intersects it; hence, ${\Deg(N)\supseteq L}$. Moreover, with the choice of coordinates of Remark \ref{possoassumere}, points of $\sigma$ are represented by skew-symmetric matrices with $(j,k)$-th entry
\[
\left\{
\begin{array}{ll}
\alpha_{i} & \mbox{if } \exists \,i \mbox{ such that } j+1=2i=k\\
-\alpha_{i} & \mbox{if } \exists \,i \mbox{ such that } k+1=2i=j\\
0 & \mbox{otherwise}
\end{array}
\right.
\]
for some sequence $(\alpha_{i})$ in $\mathbf{k}$. It is easy to check that a point of $\sigma$
\begin{itemize}
\item is a special complex if and only if $\alpha_{j}=0$ for some $j$;
\item is special of the second type if and only if $\alpha_{j}=\alpha_{k}=0$ for some $j\neq k$, i.e.~if and only if it lies on $F_{j} \cap F_{k}$.
\end{itemize}
If $N \cap F_{j} \cap F_{k}=\emptyset$ for any $j < k$, then $N$ contains $\frac{n}{2}$ special complexes of the first type and no special complexes of the second type, hence the conclusion.
\item[ii. $\Rightarrow$ i.]
Let $N$ be as in ii.. For any $i$ define $R_{i}:=N \cap \zeta^{-1}(\ell_{i})$, which is non-empty by hypothesis. We have $R_{i} \neq R_{j}$ for any $i \neq j$, as otherwise $N$ would contain a point having ${{<}}\ell_{i},\ell_{j}{{>}}$ as center. We claim that $N \subseteq \sigma$, for which it is sufficient to show that $R_{1}$ and $R_{2}$ are both contained in $\sigma$.

Let us choose the coordinates as in Remark \ref{possoassumere}. For any $i \neq j$ we have $N = {<}R_{i},R_{j}{>}$, so
\begin{equation}
\label{Runodentro}
R_{1} \cup R_{2}\subset N = \bigcap_{i < j} {<}R_{i},R_{j}{>}.
\end{equation}
Complexes in $R_{i} \subset \zeta^{-1}(\ell_{i})$ have zero $i$-th and $(i+1)$-th rows and columns, hence the entries $A_{k,l}$ of a complex $A$ in ${<}R_{i},R_{j}{>}$ are zero at least when
\[
(k,l) \in \left(\{2i-1,2i\}\times\{2j-1,2j\}\right) \cup \left(\{2j-1,2j\}\times \{2i-1,2i\}\right).
\]
From \eqref{Runodentro}, we deduce that any complex $A$ in $R_{1} \cup R_{2}$ has non-zero entries $A_{k,l}$ only if $\exists \,i$ such that $k+1=2i=l$ or $l+1=2i=k$, hence it belongs to $\sigma={{<}}H_{i}{{>}}_{1 \leq i \leq \frac{n}{2}}$.

This is enough to show that $N \subseteq \sigma$. If $N \cap F_{i} \cap F_{j}\neq\emptyset$ for some $i \neq j$, then $\Deg(N)$ would contain ${<}\ell_{i},\ell_{j}{>}$, hence a contradiction. \qedhere
\end{itemize}
\end{proof}
\end{propo}

\begin{propo}
\label{propodomcurvepari}
If $m=2$, $n\geq 4$ and $n$ is even, then $\rho$ is dominant. The general element of $\mathcal{H}$ is the union of $\frac{n}{2}$ lines spanning $\PV$. The general fiber has dimension $n-4$; its general element is a general line $N \subset \PLVs$ lying on $\sigma$, as in Proposition \ref{bohbohboh}. In particular, $\rho$ is birational if $(m,n)=(2,4)$.
\begin{proof}
We can define a rational map
\[
\xi:\xymatrix{
{\grapr(1,\PV)}^{\frac{n}{2}}
\ar@{-->}[r] &
\mathcal{H},
}
\]
defined outside the closed subset corresponding to sets of $\frac{n}{2}$ lines non-spanning $\PV$, sending $\frac{n}{2}$ lines to the corresponding point in $\mathcal{H}$. It is finite and its image $\im(\xi)$ is irreducible.

The general morphism $\phi$ has $\frac{n}{2}$ lines spanning $\PV$ as degeneracy locus, so $\im(\xi)=\im(\rho)$.
It remains to show that
\begin{equation}
\label{eqtoprove}
\dim (\im(\xi)) = \dim \mathcal{H},
\end{equation}
so that
$\overline{\im(\rho)}$ is the unique irreducible component of $\mathcal{H}$, hence $\rho$ is dominant. On the one hand, $\dim (\im(\xi))=\dim({\grapr(1,\PV)}^{\frac{n}{2}})=n^{2}-2n$. On the other, let $Y$ be the union of $\frac{n}{2}$ skew lines. For a line $\ell \subset \PV$, we have $\mathcal{N}_{\ell/\Pnm}\cong \mathcal{O}_{\mathbb{P}^{1}}^{n-2}(1)$, hence
\[
\hh^{0}(\mathcal{N}_{Y/\Pnm})=\frac{n}{2}(2(n-2))
\quad
\mbox{and} \quad \hh^{1}(\mathcal{N}_{Y/\Pnm})=0,
\]
which imply equality \eqref{eqtoprove}.

%On the other hand, let $Y$ be the union of $\frac{n}{2}$ skew lines. From
%\[
%\xymatrix{
%0 \ar[r] &
%\mathcal{T}_{Y} \ar[r] &
%\left.\mathcal{T}_{\Pnm}\right|_{Y} \ar[r] &
%\mathcal{N}_{Y/\Pnm}\ar[r] &
%0
%}
%\]
%and Euler sequence restricted to $Y$, we get
%\[
%\hh^{0}(\mathcal{N}_{Y/\Pnm})=\frac{n}{2}(2n-1-\chi(\mathcal{T}_{Y}))
%\quad
%\mbox{and} \quad \hh^{1}(\mathcal{N}_{Y/\Pnm})=0.
%\]
%Riemann-Roch Theorem yields $\chi(\mathcal{T}_{Y})=3$, so equality (\ref{eqtoprove}) holds.

%With the same technique and computations used in \cite[\textsection 4]{TanturriDegeneracy}, we get $\dim \mathcal{H}=n^{2}-2n$, which coincides with $\dim (\im(\xi))=\dim({\grapr(1,\PV)}^{\frac{n}{2}})$.

Finally, fixing a point $\cup_{i} \ell_{i}$ in $\mathcal{H}$, the general element of its preimage via $\rho$ is a general line $N \subseteq \sigma$, as showed in Proposition \ref{bohbohboh}. In particular, the space of such lines has dimension
\[
\dim \grapr(1,\sigma)=
n-4.
\]
When $n=4$, $\sigma$ is a line and the unique preimage is $\sigma$ itself.
\end{proof}
\end{propo}

We remark here that the fibers of $\rho$ can be explicitly constructed. By means of a projective transformation we can send any set of $\frac{n}{2}$ general lines to the lines chosen in Remark \ref{possoassumere}; then, we just need to apply to any line lying on $\sigma$ as in Proposition \ref{bohbohboh} the same projective transformation backwards.

%%%%%%%%%%%%%%%%%%%%

\section{\texorpdfstring{The behavior of $\rho$: odd case}{The behavior of rho: odd case}}
Let $n$ be odd from now on. The general degeneracy locus $X_{\phi}$ is easy to describe: similarly to the case $m=3$ in \cite[\textsection 6]{TanturriDegeneracy}, 
the elements in $\im(\rho)$ are the images of maps
\begin{equation}
\label{mappamappa}
\xymatrix{
\PU \ar[rr]^-{[f_{1}:\dotso:f_{n}]} && \PV,
}
\end{equation}
where $f_{1},\dotsc,f_{n}$ are forms of degree $\frac{n-1}{2}$ in $\mathbf{k}[y_{0},y_{1}]$, Pfaffians of a
general $n \times n$ skew-symmetric matrix $N$ with entries in $\mathbf{k}[y_{0},y_{1}]_{1}$. These forms are general, in the sense
that they generate the whole vector space $\mathbf{k}[y_{0},y_{1}]_{\frac{n-1}{2}}$.

\begin{lem}
\label{matrixnk}
For a general $n \times n$ skew-symmetric matrix $N$ with entries in $\mathbf{k}[y_{0},y_{1}]_{1}$, its Pfaffians of order $n-1$ span the whole $\mathbf{k}[y_{0},y_{1}]_{\frac{n-1}{2}}$.
\begin{proof}
For the $(n-1) \times (n-1)$ Pfaffians of a general $N$, not to span $\mathbf{k}[y_{0},y_{1}]_{\frac{n-1}{2}}$ is a closed condition, so it is sufficient to exhibit, for any odd $k$, a matrix $N_{k}$ not satisfying it. For this sake, we consider the $k \times k$ matrix
\[
N_{k} =
\left(
\begin{array}{ccccccc}
0 & y_{0} \\
-y_{0}& 0 & y_{1} \\
& -y_{1}& 0 & y_{0} \\
&&-y_{0}& 0 & y_{1} \\
&&&& \ddots & \ddots \\
&&&&& 0 & y_{1}\\
&&&&&-y_{1}& 0
\end{array}
\right).
\]
If we denote by $\Pf_{i}(N_{k})$ the $(k-1) \times (k-1)$ Pfaffian obtained from $N_{k}$ by deleting the $i$-th row and column, it is easy to check that
\begin{align*}
& \Pf_{2i+1}(N_{k})= y_{0}^{i} y_{1}^{\frac{k-1}{2}-i} & \mbox{for any } 0 \leq i \leq \frac{k-1}{2},\\
& \Pf_{2i}(N_{k})=0 & \mbox{for any } 1 \leq i \leq \frac{k-1}{2},
\end{align*}
and this concludes the proof.
\end{proof}
\end{lem}

\begin{rmk}
\label{tuttipfaff}
As a consequence of the previous lemma, every sequence of general forms $f_{1},\dotsc,f_{n}$ of degree $\frac{n-1}{2}$ corresponds to the sequence of Pfaffians of a suitable skew-symmetric matrix. Indeed, these forms can be expressed as linear combination of the Pfaffians of $N_{k}$ above, giving rise to $\beta \in \PGL(V)$ such that the diagram
\[
\xymatrix{
& & Y_{1} \ar^-{\beta}[drr]\\
\PU \ar[urr]^-{[f_{1}:\dotso:f_{n}]\,}_-{\sim} \ar[rrrr]_-{[\Pf_{1}(N_{k}):\dotso:\Pf_{k}(N_{k})]\,}^-{\sim}
& & & & Y_{2}
}
\]
commutes. This produces an automorphism of $\PU$, hence a change of basis of $\mathbf{k}[y_{0},y_{1}]_{1}$. In terms of this new basis, $N_{k}$ has the desired Pfaffians $f_{1},\dotsc,f_{n}$.
\end{rmk}

\begin{propo}
If $m=2$, $n \geq 5$ and $n$ is odd, then $\rho$ is dominant. The general element of $\mathcal{H}$ is the image in $\PV$ of a map
\[
\xymatrix{
\PU \ar[rr]^-{[f_{1}:\dotso:f_{n}]} && \PV,
}
\]
where $f_{1},\dotsc,f_{n}$ are forms of degree $\frac{n-1}{2}$ spanning $\mathbf{k}[y_{0},y_{1}]_{\frac{n-1}{2}}$. The general fiber of $\rho$ has dimension $\frac{n^{2}-3n}{2}$.
\begin{proof}
Let $r=\dim(\mathbf{k}[y_{0},y_{1}]_{\frac{n-1}{2}})-1=\frac{n-1}{2}$. We can define a rational map
\[
\xi:
\xymatrix{
\mathbb{A}^{(r+1)n} \ar@{-->}[r] & \mathcal{H}
}
\]
sending
an $n$-tuple of forms $f_{1},\dotsc,f_{n}$ to the image of the map \eqref{mappamappa}. It is defined on the $n$-tuples which span the whole linear space $\mathbf{k}[y_{0},y_{1}]_{\frac{n-1}{2}}$; its image $\im(\xi)$ is irreducible and its dimension is easily computable. Indeed, on the one hand there is a natural $\GL_{2}$-action on $\mathbf{k}[y_{0},y_{1}]_{1}$, acting as a change of basis on $U$; this induces an action on $\mathbf{k}[y_{0},y_{1}]_{\frac{n-1}{2}}$ and therefore on $\mathbb{A}^{(r+1)n}$, and one can see that $\xi$ 
factors through this action. On the other, take two points $Y_{1},Y_{2}$ in $\im(\xi)$ such that $Y_{1}=Y_{2}$.
By the commutativity of the diagram
\[
\xymatrix{
& & Y_{1} \ar@{=}[dd]\\
\PU \ar[urr]^-{[f_{1}:\dotso:f_{n}]\,}_-{\sim} \ar[drr]_-{[g_{1}:\dotso:g_{n}]\,}^-{\sim}\\
& & Y_{2}
}
\]
we get an automorphism of $\PU$, i.e.~the two maps $[f_{1}:\dotso:f_{n}]$ and $[g_{1}:\dotso:g_{n}]$ belong to the
same class modulo $\GL_{2}$. Hence
%\begin{align*}
%\dim (\im(\xi))	&	=\dim (\mathbb{A}^{(r+1)n}) - \dim (\GL_{2}) \\
%				&	= \frac{n^{2}+n-8}{2}. 
%\end{align*}
\begin{equation*}
\dim (\im(\xi))	=\dim (\mathbb{A}^{(r+1)n}) - \dim (\GL_{2}) = \frac{n^{2}+n-8}{2}. 
\end{equation*}

Since $n$ general forms are obtained as the Pfaffians of a suitable matrix $N$ by Remark \ref{tuttipfaff}, $\im(\rho)=\im(\xi)$. If we prove that
\begin{equation}
\label{eqtoprove2}
\dim (\im(\xi)) = \dim \mathcal{H},
\end{equation}
then $\overline{\im(\rho)}$ is the unique irreducible component of $\mathcal{H}$, hence $\rho$ is dominant. Let $Y \in \im(\xi)$. From
\[
\xymatrix{
0 \ar[r] &
\mathcal{T}_{Y} \ar[r] &
\left.\mathcal{T}_{\Pnm}\right|_{Y} \ar[r] &
\mathcal{N}_{Y/\Pnm}\ar[r] &
0
}
\]
and Euler sequence restricted to $Y$, we get
\[
\hh^{0}(\mathcal{N}_{Y/\Pnm})=n\frac{n+1}{2}-1-\chi(\mathcal{T}_{Y})
\quad
\mbox{and} \quad \hh^{1}(\mathcal{N}_{Y/\Pnm})=0.
\]
Riemann-Roch Theorem yields $\chi(\mathcal{T}_{Y})=3$, so equality \eqref{eqtoprove2} holds.

%This can be done again as in the proof of Proposition \ref{propodomcurvepari}.

The dimension of the fibers is finally
\[
\dim \GRAd - \dim \mathcal{H} = \frac{n^{2}-3n}{2}. \qedhere
\]

\end{proof}
\end{propo}

As in the even case, also for odd values of $n$ it is possible to construct explicitly the fibers of $\rho$. Up to a projective transformation, we can assume that the degeneracy locus is the image of the Pfaffians of the matrix $N_{k}$ from Lemma \ref{matrixnk}. We can find by linear algebra all the skew-symmetric matrices having those as Pfaffians and apply the same projective transformation backwards to get the elements of the desired fiber.


\begin{thebibliography}{Tan15}

\bibitem[BM01]{BazanMezzetti}
D.~Bazan and E.~Mezzetti.
\newblock On the construction of some {B}uchsbaum varieties and the {H}ilbert
  scheme of elliptic scrolls in {$\Bbb P^5$}.
\newblock {\em Geom. Dedicata}, 86(1-3):191--204, 2001.

\bibitem[FF10]{FaenziFania}
D.~Faenzi and M.~L.~Fania.
\newblock Skew-symmetric matrices and {P}alatini scrolls.
\newblock {\em Math. Ann.}, 347(4):859--883, 2010.

\bibitem[FM02]{FaniaMezzetti}
M.~L.~Fania and E.~Mezzetti.
\newblock On the {H}ilbert scheme of {P}alatini threefolds.
\newblock {\em Adv. Geom.}, 2(4):371--389, 2002.

\bibitem[Ott92]{Ottaviani}
G.~Ottaviani.
\newblock On {$3$}-folds in {$\bold P^5$} which are scrolls.
\newblock {\em Ann. Scuola Norm. Sup. Pisa Cl. Sci. (4)}, 19(3):451--471, 1992.

\bibitem[Tan15]{TanturriDegeneracy}
F.~Tanturri.
\newblock On the {H}ilbert scheme of degeneracy loci of twisted differential
  forms.
\newblock To appear in {\em Trans. Amer. Math. Soc.}.

\end{thebibliography}
\end{document}